\documentclass{tac}
\usepackage{xy}
\xyoption{all}

\def\mathcaldef#1{\expandafter\def\csname#1\endcsname{{\cal#1}}}

\newcommand{\tto}{\mathop{\to}\limits}
\newcommand{\mapstto}{\mathop{\longmapsto}\limits}
\def\too{\longrightarrow}
\def\mapstoo{\longmapsto}
\def\q{\quad}
\def\qq{\quad\quad}
\def\qv{\qq ;\qq}
\def\v{``}
\def\imp{\,\Rightarrow\,} 
\def\iff{\q\Longleftrightarrow\q}
\def\iso{\,\cong\,}
\def\equ{\,\simeq\,}
\def\inc{\hookrightarrow}
\def\down{\downarrow\!\!}
\def\up{\uparrow\!\!}
\def\la{\langle}
\def\ra{\rangle}
\def\adj{\dashv}
\def\op{^{\rm op}}
\def\ort{^{\perp}}
\def\ex{\exists}
\def\comp{\Gamma_{!}}

\def\lam{\lambda}
\def\De{\Delta}
\def\de{\delta}
\def\ne{\nu}
\def\eps{\varepsilon}
\def\ov{\overline}

\mathcaldef{C}
\mathcaldef{E}
\mathcaldef{M}
\mathcaldef{S}
\mathcaldef{P}
\mathcaldef{A}

\mathbfdef{Set}
\mathbfdef{Cat}
\mathbfdef{Gph}
\mathbfdef{Pos}
\mathbfdef{y}

\mathrmdef{hom}
\mathrmdef{id}
\mathrmdef{Nat}
\mathrmdef{Lim}
\mathrmdef{Colim}
\mathrmdef{colim}
\mathrmdef{Cone}
\mathrmdef{tot}

\def\CatX{\Cat\!/\! X}

\def\PrX{\Set^{X\op}}

\def\CX{\C\!/\! X}
\def\CY{\C\!/Y}

\def\MX{\M\!/\! X}
\def\MY{\M\!/Y}

\def\MPX{\M\!/\P X}

\def\Mu{\M\!/1}
\def\MoX{\M'\!/\! X}
\def\MoY{\M'\!/Y}
\def\EM{{\cal(E,M)}}

\newtheorem{prop}{Proposition}
\newtheorem{corol}{Corollary}
\let\pf\proof
\let\epf\endproof
\def\eq{\begin{equation}}
\def\eeq{\end{equation}}
\def\eqa*{\begin{eqnarray*}}
\def\eeqa*{\end{eqnarray*}}

\author{Claudio Pisani}
\address{via Gioberti 86,\\ 10128 Torino, Italy.}
\title{Categories of categories}
\copyrightyear{2007}
\eaddress{pisclau@yahoo.it}

\begin{document}

\maketitle

\begin{abstract}
A certain amount of category theory is developed in an arbitrary finitely complete category 
with a factorization system on it, playing the role of the comprehensive factorization system on $\Cat$.
Those aspects related to the concepts of finality (in particular terminal objects), discreteness and components,
representability, colimits and universal arrows, seem to be best expressed in this very general setting.
Furthermore, at this level we are in fact doing not only $\EM$-category theory but, in a sense, 
also $\EM$-topology. 
Other axioms, regarding power objects, duality, exponentials and the arrow object, are considered.   
\end{abstract}

\tableofcontents

\section{Introduction}
\label{intro}

The basic logic of a factorization system $\EM$ on a finitely complete category $\C$ 
has been known for a long time:
it gives rise to a (variant of) elementary existential doctrine with the comprehension scheme
(see \cite{law70} and~\cite{kor77}); or, otherwise stated, to a subfibration of the codomain fibration 
on $\C$, which is itself a bifibration.
Indeed, substitution $\De_f:\MY\to\MX$ along $f:X\to Y$ is given by pullbacks, and since
the category $\MX$ of predicates over $X\in\C$ is a reflective full subcategory of $\CX$:
\[ \down \,\,\,\adj i : \MX \inc \CX \]
the existential quantification $\ex_f \adj \De_f$ can be obtained as $\ex_f m = \,\,\down (f\circ im)$.

In particular, the doctrine associated to the comprehensive factorization system 
\[ \MX \equ \PrX \qq X\in\Cat \]  
(where truth values are sets and existential quantification corresponds to Kan extensions) was among the 
instances of hyperdoctrine in~\cite{law69} and~\cite{law70}.
The comprehensive factorization system itself was made explicit in~\cite{stw73}. 
(We here consider the dual choice: $\E$ and $\M$ are the classes of final functors and discrete fibrations,
rather than initial functors and discrete opfibrations).

About at the same time, the relevance of the \v comprehension adjunction" 
\[ \down \,\,\,\adj i : \MX \equ\PrX \inc \CatX \]
in the calculus of colimits was shown thoroughly in~\cite{par73}.
By developing further these items in recent years (see~\cite{pis}), it became clear that  
several categorical properties rest on the above reflection, that is on the comprehensive factorization
system. In fact, in the present paper we will show how some of these properties, when properly formulated,  
{\em depend only on the general logic of factorization systems}.

So, category theory has a kernel which can be developed in any $\EM$-category $\C$, that is in
any finitely complete category with a factorization system on it.
However, in general the objects of $\C$ do not yet deserve to be called categories: 
in order to consider $\C$ a \v category of categories" we should assume at least an \v arrow object"
as in Section~\ref{ar} (see also below).

On the contrary, other axioms on $\C$ (see Section~\ref{top}) would lead to consider it
rather as a category of \v topological spaces", including some kind of \v infinitesimal" objects.
(The term \v space" is in some respects inadequate: for instance, although a right adjoint to the 
inclusion $\MX\inc\CX$ does exist for $\C = \Cat$, we presently are not assuming it,
so that $\C$ is not a category of cohesion relative $\S = \Mu$ in the sense of~\cite{law07}.)

In such a \v topological" category, any space over a base space $X$ has a {\em discrete} (or {\em \'etal\'e})
reflection in $\MX$. In particular, any space has a \v set of components", 
and the discrete reflection of a point $x:1\to X$ is the (infinitesimal) 
{\em neighborhood} $\ne(x):N\to X$ of $x$.
The full subcategory $\ov X\inc\MX$ generated by the neighboroods of the points of $X$ is the
{\em adherence} category of $X$. 
If a space $p:P\to X$ over $X$ is in a neighborhood $\ne(x)$,
then \v $p$ is near the point $x$"; more generally, if neighbhorhoods are not monomorphisms,
a map $p\to\ne(x)$ \v is a way to be near the point $x$".
For a discrete space $p:P\to X$ in $\MX$, the following property is now almost tautological:
for any $x\in X$ and any $a\in P$ over it, the neighborhood $\ne(x)$ can be uniquely lifted 
to the neighborhood $\ne(a)$ of $a$; equivalently, for any $a\in P$, $f\circ\ne(a) = \ne(pa)$.
Thus $p:P\to X$ is indeed a \v local hemeomorphism" in a very natural sense.
(The fact that local homeomorphisms are the topological correspective of the discrete fibrations
in $\Cat$ is further supported by~\cite{hof05}, where it is shown that \'etal\'e spaces
coincide with pullback-stable \v discrete ultrafilter fibrations"; 
in our context, stability of discreteness always holds.)
 
In $\Cat$, the neighborhood (in the above sense) of the point-object $x\in X$ 
is the discrete fibration $X/x$ corresponding to the presheaf represented by $x$, 
and a way for a functor $p:P\to X$ to be near $x$ is simply a cone of base $p$ and vertex $x$.
So, a reflection of a functor $p:P\to X$ in $\CatX$ in the adherence category $\ov X\iso X$ is a colimit of $p$. 
Furthermore, $f:X\to Y$ has a right adjoint iff the inverse image (pullback) under $f$ of any neighborhood in $Y$ 
is a neighborhood in $X$. 

In Section~\ref{fs}, we briefly review factorization systems in general, and the 
comprehensive factorization system on $\Cat$.
In Sections~\ref{gen} and~\ref{thm} we define and prove some classical concepts and
theorems of category theory in the general context of $\EM$-categories.

In the rest of the paper we introduce some axioms, mostly aiming at a more accurate description
of the category of categories. 
In Section~\ref{pow}, we define power objects (playing the role of presheaf categories) and prove that
if $X$ has such a \v Yoneda map" $\y:X\to\P X$, the classical definition-characterization of finality holds:
a map $e:P\to X$ is in $\E$ iff for any map $f:X\to Y$ the colimits of $f$ and $f\circ e$
are the same (either existing if the other one does).
If the power object $\Omega := \P 1$ exists in $\C$ it may be considered the \v internal sets" 
or \v truth values" object.
In Section~\ref{dual}, we assume a duality functor $(-)':\C\to\C$ and exponentials, 
so that we have \v two-sided" category theory. In particular \v hom maps" $X'\times X\to\Omega$ 
can be defined as those maps such that both the transposes are Yoneda maps.

In Section~\ref{ar}, we analyze the case where $\EM$ is generated by a pointed object $t:1\to T$, 
so that the discrete spaces $m:M\to X$ over $X$ coincide with the \v discrete $T$-fibrations":
for any $T$-figure $f:T\to X$ of the base $X$ and any point $a$ of $M$ over $ft$,
there is a unique lifting $f':T\to M$ of $f$ with $f't = a$.
If $T$ is in fact bipointed $s,t:1\to 2$, and $s$ generates the dual of $\EM$, we call it an \v arrow object".
Supposing that Lawvere's axioms on the pushouts $3$ and $4$ of $2$ hold, any $\C$-space $X\in\C$ 
is a category $X^\star$ (see~\cite{law66}), and there is a natural transformation $X^\star\to\ov X:\C\to\Cat$. 

Most of the classical examples of factorization systems give rather uninteresting 
instances of \v category theories". (Indeed, if all the points $1\to X$ are in $\M$, 
the adherence categories are discrete and only the constant maps can have (co)limits). 
However, apart from categories and \v topological spaces", there are at least two $\EM$-categories which may take
advantage from being considered in the present perspective: the category of (reflexive) graphs
and the category of posets, with $\M$ given respectively by graph fibrations and by
inclusions of lower-sets.
These are briefly considered in Section~\ref{ex},
where we also discuss some properties of \v topological categories".

Some of these topics have been presented at the International Conference on Category Theory (CT2007) 
held at Carvoeiro in June 2007.

\section{The comprehensive factorization system}
\label{fs}

We review the basic facts regarding factorization systems in general and the 
comprehensive factorization system on $\Cat$, on which our theory is modeled.

\subsection{The logic of a factorization system}

We assume throughout that $\C$ is a finitely complete category.
Recall that the {\bf orthogonality} relation between the arrows of $\C$ 
is defined as follows: if $e:E\to X$ and $m:M\to Y$, then $e\perp m$ iff
every commutative square $ge = mf$ has a unique diagonal:
\eq  \label{1fs}
\xymatrix@R=3pc@C=3pc{
E \ar[r]^f\ar[d]_e          &  M \ar[d]^m \\
X \ar[r]^g\ar@{..>}[ur]^u   &  Y           }
\eeq
Thus we have a Galois connection on the classes of arrows of $\C$.
If $\E$ and $\M$ are two such classes, stable with respect to $\perp$, 
\[ 
\E = \,^\perp\!\M \qv \M = \E^\perp 
\]
one says that the pair $\EM$ is a {\bf pre-factorization system} on $X$;
if furthermore any $f$ in $\C$ factors as $f = me$ with $e\in\E$ and $m\in\M$,
then $\EM$ is a {\bf factorization system}.
Note that for any class $\A$ of arrows in $\C$, $(^\perp\!(\A^\perp),\A^\perp)$
is a pre-factorization system on $X$: we say that it is \v generated by $\A$".
The following properties are standard:
\begin{prop}  \label{fs1}
Let $\EM$ be a pre-factorization system on $\C$.
\begin{enumerate}
\item
An arrow $f$ is both in $\E$ and in $\M$ iff it is an isomorphism.
\item
If $m$ and $m'$ are consecutive arrows with $m\in\M$, then $m'\in\M$ iff $m'\circ m\in\M$.
\item
If $m\in\M$, the pullback $f^*m$ along any map is also in $\M$.
\end{enumerate}
\end{prop} 
\epf
\noindent Of course, the dual properties concerning $\E$ hold as well.
\begin{corol}  \label{fs2}
If $\EM$ is a pre-factorization system on $\C$, then
\begin{enumerate}
\item
$\E$ and $\M$ are (lluf) subcategories of $X$.
\item
The inclusions \( i_X:\MX \inc \CX \) are full.
\item
By restricting the codomain fibration to the arrows in $\M$, one obtains a subfibration.
\end{enumerate}
\end{corol}  \label{fs3}
\epf
\noindent So, for any $f:X\to Y$ in $\C$, the pullback functor $f^*: \CY\to\CX$ 
restricts to
\eq
\De_f: \MY\to\MX
\eeq
\begin{prop}   \label{fs4}
Let $\EM$ be a pre-factorization system on $\C$ and $p:P\to X$.
The following are equivalent:
\begin{enumerate}
\item
$p = m\circ e$ with $m\in\M$ and $e\in\E$. 
\item
The map $e:p\to m$ over $X$ is a reflection of $p\in\CX$ in the subcetegory $\MX$. 
\end{enumerate}
\end{prop}
\epf
\begin{corol}   \label{fs5}
A pre-factorization system $\EM$ on $\C$ is a factorization system
if and only if $\MX$ is reflective in $\CX$ for any $X\in\C$.
\end{corol}
\epf
Thus, the key facts about factorization systems can be summarized as follows
\begin{itemize}
\item
A factorization $p = m\circ e$, with $m\in\M$ and $e\in\E$,                    
gives the following universal property: 
any map in $\CX$ from $p$ to $n$, with $n\in\M$, factors uniquely through $e$: 

\eq    \label{2fs}
\xymatrix@R=4pc@C=4pc{
P \ar[r]^e\ar[dr]_p\ar@/^1.5pc/[rr]  & M \ar[d]^m \ar@{..>}[r] & N \ar[dl]^n \\
                                     & X                       &              }
\eeq
So if a factorization $p = me$ of any map $p$ has been fixed, we have reflections 
\[ \down_X\adj i_X : \MX\to\CX \]
with $\down p = m$; conversely, such a reflection yields factorizations $p = \,\down p\circ e$.
\item
In particular, an arrow is in $\E$ iff its reflection is terminal in $\MX$ 
(that is, iff it is an isomorphism): 
\eq   \label{3fs}
e\in\E \iff \down e \iso 1_X.
\eeq
\item
The fibration $\MX$ is a subfibration of the codomain fibration $\CX$, 
\eq   \label{4fs}
\xymatrix@R=4pc@C=4pc{
\MX  \ar[d]^i  & \MY \ar[l]_{\De_f}\ar[d]^i   \\
\CX            & \CY \ar[l]_{f^*}              }
\eeq
and is itself a bifibration
\[ \ex_f\adj\De_f \]
where the existential quantification can be obtained as 
\[ \ex_f m \iso \down f_!(im) \]
Since the left adjoint square also commutes (up to isomorphisms) 
\eq   \label{4bfs}
\xymatrix@R=4pc@C=4pc{
\MX \ar[r]^{\ex_f}           & \MY                  \\
\CX \ar[u]_\down\ar[r]^{f_!}    & \CY  \ar[u]_\down   }
\eeq  
we get 
\eq   \label{5fs}
\down(f\circ p) \iso \ex_f\down p
\eeq
for any $p\in\CX$. 
(Recall that $f_! \adj f^*$ is given by composition with $f$.)
\item
In particular, if $e:P\to X$ is in $\E$ then for any $f:X\to Y$
\eq   \label{6fs}
\down(f\circ e) \iso \ex_f 1_P \iso \,\down f
\eeq
So, $e:P\to X$ is in $\E$ iff, for any $f:X\to Y$, $f$ and $fe$ have the same reflection.
\end{itemize}

\subsection{The comprehensive factorization system}

Let the functor $t:1\to 2$ be the codomain of the arrow category,
and let $\EM$ be the pre-factorization system generated by $t$. 
Then the class $\M = \,\ort t$ of functors orthogonal to $t$, is the (lluf) subcategory of 
discrete fibrations, and the comprehensive factorization system reduces to the following facts:
\begin{itemize}
\item
for any category $X$, $\MX$ is reflective in $\CatX$;
\item
the class $\E = \M\ort$ is the (lluf) subcategory of final functors.
\end{itemize}
Both of them are a consequence of the following proposition, whose proof is an easy generalization 
of the Yoneda Lemma (see~\cite{pis2}): 
\begin{prop}   \label{fs6}
The left adjoint $\down\,\,\,\adj i$ to the full inclusion $i:\MX \inc \CatX$
is given on a $p:P\to X$ by (the discrete fibration corresponding to the presheaf)
\eq    \label{7fs}
\down p \iso \comp(-/p) 
\eeq
\end{prop}
\epf
In fact, when $P = 1$ we get the reflection of an object $x:1\to X$,
\[ \down x = X/x \]
and the universal property of the reflection reduces to the 
(discrete fibration version of the) Yoneda Lemma; 
while when $X = 1$ the reflection
\[ \down\,\,\, : \Cat \to \Mu \equ \Set \]
is simply the component functor $\comp$.

Now, by~(\ref{3fs}), a functor $p:P\to X$ is in $\E$ iff $\comp(x/p) = 1$ for any $x\in X$,
and this is the definition of final functors given in~\cite{mac71}). 
The other classical characterization of finality via colimits 
will be proved in a more general context in Section~\ref{pow}.

\section{The general context}
\label{gen}

In this section we define several categorical concepts in the context of arbitrary $\EM$-categories.
Thus, we assume that $\C$ is a finitely complete category with a factorization system on it.
We also suppose that a reflection $\,\,\down_X:\CX\to\MX$, for any $X\in\C$ 
(or, equivalently a factorization for any arrow in $\C$) has been fixed.
The objects and the arrows of $\C$ will be called {\bf $\C$-spaces} (or simply spaces) and {\bf maps}.

\subsection{Final maps, discrete spaces over a space and sets}

A map in $\E$ is {\bf final}.
An object $m:M\to X$ of $\MX$ is a {\bf discrete space} over $X$.
A {\bf $\C$-set} is a discrete space over $1$, and we denote by $\S := \Mu\inc\C$  
the reflective full subcategory of $\C$-sets, which are the truth values of our doctrine.
(We will often drop the prefix $\C$.)

\subsection{Components}

The reflection $\comp:=\,\,\down_1:\C\to\S$ is the {\bf components} functor. 
So the factorization of the terminal map $!_X$
\[
\xymatrix@R=4pc@C=4pc{
X  \ar[r]^e\ar[dr]_{!}  & S \ar[d]^m \\
                        & 1           }
\]
gives the set $S = \comp X$ of components of $X$, with the universal property:

\eq   \label{1gen}
\xymatrix@R=4pc@C=4pc{
X \ar[r]^e\ar[dr]\ar@/^1.5pc/[rr]  & S \ar[d]^m \ar@{..>}[r] & S' \ar[dl]^n \\
                                   & 1                       &               }
\eeq
In particular, a space $X$ is {\bf connected} if $\comp X = 1$ (or, equivalently, if the 
terminal map is final). In this case we have the familiar universal property:

\eq   \label{2gen}
\xymatrix@R=4pc@C=4pc{
X \ar[r]^e\ar[dr]^e\ar@/^1.5pc/[rr]  & 1 \ar[d] \ar@{..>}[r] & S \ar[dl] \\
                                   & 1                     &               }
\eeq
that is, any map to a set is constant.
Note that by Proposition~\ref{fs1} a space with a final point is connected.

\subsection{Neighborhoods}
\label{ne}

Given a point $x:1\to X$, the {\bf neighborhood} of $x$ is the reflection $\ne(x) :=\,\down x$
of $x$ in $\MX$. So the factorization of the point $x$ 
\[
\xymatrix@R=4pc@C=4pc{
1  \ar[r]^e\ar[dr]_x    & N \ar[d]^{\ne(x)} \\
                        & X           }
\]
gives the neighborhood $m = \ne(x)$ of $x$ with a final point of its total.
Thus, a discrete space $m:M\to X$ over $X$ is isomorphic to a neighborhood 
iff its total $M$ has a final point $e:1\to M$, and in this case $m \iso \ne(me)$.
(In particular, any $\C$-space is \v locally connected".)

The universal property of discrete reflection now becomes:
\[
\xymatrix@R=4pc@C=4pc{
1 \ar[r]^e\ar[dr]_x\ar@/^1.5pc/[rr]^a  & N \ar[d]^{\ne(x)} \ar@{..>}[r]^u & M \ar[dl]^m \\
                                       & X                                &              }
\]
That is, given any point $a$ over $x$  of the total of a discrete space over $X$,
there is a unique map $\ne(x)\to M$ over $X$ which takes the final point of its total to $a$.
Note that $u$ is itself discrete, so that it is the neighborhood $u = \ne(a)$ of the point $a$ of $M$. 
(This is the \v local homeomorphism property" mentioned in the Introduction.)
Otherwise posed, we may rearrange the diagram as follows:
\[
\xymatrix@R=3.5pc@C=3.5pc{
1 \ar[r]^a\ar[d]_e                &  M \ar[d]^m   \\
N \ar[r]^{\ne(x)}\ar@{..>}[ur]^u  &  X             }
\]
so that $\ne(x)$ is characterized by the fact that it is {\em initial} among the {\em discrete} spaces
over $X$ with a (distinguished) point over $x$.
The dual characterization is also significant:
the neighborhood of $x$ is {\em final} among the spaces over $X$ with a (distinguished) {\em final} point over $x$:
\[
\xymatrix@R=3.5pc@C=3.5pc{
1 \ar[r]^e\ar[d]_t         &  N \ar[d]^{\ne(x)} \\
T \ar[r]^f\ar@{..>}[ur]^u  &  X                    }
\]
Intuitively, a $\C$-space $T$ with a final point $t:1\to T$ is \v absolutely concentrated" around $t$, 
so that a map $f:T\to X$ parameterizes a portion of the space $X$ concentrated around $ft\in X$. 
The neighborhood $\ne(x):N\to X$ is thus characterized as the biggest of these portions. 
In fact, we will see that in the above diagram $x$ is the absolute colimit of $f$ in a precise sense 
(which becomes the usual one in $\Cat$), with $u$ as colimiting cone.

\subsection{Neighborhoods in $\Cat$}

As mentioned before, in $\Cat$ the neighborhood of the point-object $x$ 
is the discrete fibration $\ne(x) : X/x \to X$ (with the terminal object $\id_x:1\to X/x$ as reflection map)
corresponding to the presheaf represented by $x$.
Thus, a category $p:P\to X$ over $X$ is isomorphic (over $X$) to $\ne(x):X/x\to X$ iff it is a 
discrete fibration and $P$ has a terminal object $e:1\to P$ over $x$. 
The isomorphism is unique if it is required to send $e$ to $\id_x$.  
The universal property of neighborhoods now becomes  
(the discrete fibration version of) the Yoneda Lemma:
a discrete fibration $n:N\to X$ over $X$ with an object $e:1\to N$ over $x\in X$
is isomorphic to $X/x \to X$ (with the object $\id_x$ over $x$) iff
for any discrete fibration $m:M\to X$, composition with $e$ gives a bijection between 
the functors $N\to M$ over $X$ and the objects $a$ of $M$ over $x$. 
\eq   \label{3gen}
\xymatrix@R=3.5pc@C=3.5pc{
1 \ar[r]^e\ar[dr]_x\ar@/^1.5pc/[rr]^a  & N \ar[d]^n \ar@{..>}[r] & M \ar[dl]^m \\
                                       & X                       &              }
\eeq
(Of course, $e$ is then a terminal object and is usually known as a \v universal element" 
of the presheaf corresponding to $n$.)

The dual characterization is less widely known:
a category $n:N\to X$ over $X$ with a final object $e:1\to N$ over $x\in X$
is isomorphic to $X/x \to X$ (with the object $\id_x$ over $x$) iff 
for any category with a final object $t:1\to T$, composition with $n$ gives a bijection 
between the functors $T\to N$ sending $t$ to $e$ and the functors $f:T\to X$ sending $t$ to $x$. 
\eq   \label{3bgen}
\xymatrix@R=3.5pc@C=3.5pc{
1 \ar[r]^e\ar[dr]_x\ar@/^1.5pc/[rr]^t  & N \ar[d]^n & T \ar@{..>}[l]\ar[dl]^f \\
                                       & X          &                         }
\eeq

The two dual characterizations of neighborhoods become evident in the $\EM$-category  
of posets (see Section~\ref{pos}). 
Given an element $x\in X$ of a poset, the principal sieve $\down x$ is both 
the smallest lower-set containing $x$ and the biggest part of $X$ having $x$ has a maximum.

\subsection{Pullbacks, fibers and constant spaces}

In the following pullbacks,
\eq   \label{4gen}
\xymatrix@R=4pc@C=4pc{
Mx \ar[r]\ar[d]_{x^*m}  & M \ar[d]^m \\
1  \ar[r]^x             & X          }
\qq\qq            
\xymatrix@R=4pc@C=4pc{
M  \ar[r]\ar[d]_{!^*m}  & S \ar[d]^m \\
X  \ar[r]^{!}           & 1          }
\eeq
the set $Mx$ is the (discrete) {\bf fibre} of $m$ over $x$ (or also the {\bf value} of $m$ at $x$),
while $!^*m$ is the {\bf constant} (discrete) space $S\times X \to X$ over $X$ (with constant value $S$).

\subsection{The adherence pseudo-functor}

The {\bf adherence category} $\ov X$ of $X$ is the full subcategory of $\MX$
generated by the neighborhoods of the points of $X$.
We will often denote the object $\ne(x)$ of $\ov X$ simply by $x$. 

Note that in $\Cat$, $\ov X \iso X$.
Note also that, by Proposition~\ref{fs1}, any point of a $\C$-set $S$ is discrete over $S$, 
so that $\ov S$ is a $\Cat$-set (a discrete category). 

By~(\ref{5fs}), we get in particular that for any map $f:X\to Y$
\[  
\ex_f \down x \iso\, \down(f_!x) 
\]
Then the direct image (that is, the existential quantification, not to be confused with $f_!$) 
along $f$ preserves neighborhoods:
\eq   \label{5gen}
\ex_f \ne(x) \iso \ne(fx) 
\eeq
so that $\ex_f$ restricts to a functor $\ov f:\ov X\to\ov Y$, and we have a pseudo-functor
\eq   \label{6gen}
\ov{(-)}:\C\to\Cat 
\eeq

\subsection{Displacements and cones}
\label{disp}

Given two points $x$ and $y$ in $X$, a {\bf displacement} $\lam$ from $x$ to $y$ is a map $x\to \ne(y)$ 
over $X$ from $x$ to the neighborhood of $y$:  
\eq   \label{7gen}
\xymatrix@R=4pc@C=4pc{
1  \ar[r]^\lam\ar[dr]_x    & N \ar[d]^{\ne(y)} & 1 \ar[l]_e\ar[dl]^y \\
                           & X                 &                        }
\eeq
Note that the displacements from $x$ to $y$ correspond to the arrows $\ne(x) \to \ne(y)$ in $\ov X$, 
so that composition in $\ov X$ can be seen as a Kleisli construction.

More generally, given a space $p:P\to X$ over $X$ and a point $y$ in $X$, a {\bf cone} 
with base $p$ and vertex $y$ is a map $\lam : p\to \ne(y)$ over $X$ from $p$ to the neighborhood of $y$:
\eq   \label{8gen} 
\xymatrix@R=4pc@C=4pc{
P  \ar[r]^\lam\ar[dr]_p    & N \ar[d]^{\ne(y)} & 1 \ar[l]_e\ar[dl]^y \\
                           & X                 &                        }
\eeq
It is easy to see that in $\Cat$ we get the usual cones (see~\cite{pis}).
In general, a cone expresses a way for its base to be {\em near} the vertex.

\subsection{Colimits}

A cone $\lam:p\to \ne(x)$ is {\bf colimiting} if it is a reflection of $p$
in the adherence category $\ov X$.
In this case, we say that the vertex $x\in\ov X$ is the {\bf colimit} of $p:P\to X$.
Intuitively, a colimiting cone expresses the \v best way" for $p$ to be near a point.
Of course, in $\Cat$ we obtain the usual concept of colimit of a functor.

We in fact define two colimit functors 
\eq   \label{9gen}
\Colim_X:\CX\to\ov X \qv \colim_X:\MX\to\ov X 
\eeq
as the partially defined left adjoints to the full inclusions
\[
k_X:\ov X\to\CX \qv i_X:\ov X\to\MX 
\]
Since the diagrams
\[ 
\xymatrix@1@C=3.5pc{\ov X\ar[r]^j \ar@/_1.5pc/[rr]_k & \MX \ar[r]^i & \CX  } 
\qv
\xymatrix@1@C=3.5pc{\ov X & \MX \ar[l]_{\colim} & \CX \ar[l]_\down\ar@/^1.5pc/[ll]^{\Colim} } 
\]
commute (up to isomorphisms), the Colimit of a map $p\in\CX$ is determined by its reflection:
\eq   \label{10gen}
\Colim\, p \iso \colim \down p 
\eeq
In fact, any cone $\lam : p\to \ne(x)$ has a {\bf kernel cone} $\lam' : \,\down p\to \ne(x)$
with a discrete domain; the cone is universal iff its kernel is so: 
\eq   \label{11gen} 
\xymatrix@R=2.5pc@C=4pc{
p \ar[r] \ar[dr]_\lam & \down p \ar@{..>}[d]^{\lam'}  \\
                      & \ne(x)                         }
\eeq 
As a consequence of~(\ref{5fs}) and~(\ref{5gen}) we also have
\eq   \label{10bgen}
\Colim\,(f\circ p) \iso \colim\, \ex_f\down p 
\eeq

\subsection{Preservation of colimits}

Given a cone with a discrete base $\lam : m\to \ne(x)$ and a map $f:X\to Y$,
we have by~(\ref{5gen}) an {\bf image cone} under $f$:
\eq   \label{12gen}  
\ex_f\lam : \ex_f m \to \ne(fx) 
\eeq
We say that $f$ {\bf preserves the colimiting cone} $\lam : m\to \ne(x)$ 
if the image cone $\ex_f\lam$ is also colimiting.
A map $f:X\to Y$ is {\bf colimit preserving} if it preserves all the limiting cones with a discrete base.

In $\Cat$, a functor preserves a colimiting cone in the classical sense 
iff it preseves its kernel cone in the above sense: 
\begin{prop}
Let $p:P\to X$ and $f:X\to Y$ be functors, and $\lam:p\to\ne(x_0)$ a cone.
Then the kernel cone of the \v classical" image cone of $\lam$ along $f$  
is the direct image of the kernel of $\lam$.
\end{prop}
\pf
The \v classical" image of $\lam_a:pa\to x_0,\,a\in P$ is of course the cone 
\[ f\lam_a : fpa\to fx_0 ,\,\,a\in P \]
Straightforward calculations (see also~\cite{pis2}) show that: 
\begin{itemize}
\item
The kernel cone $\lam':\,\down p\to\ne(x)$ of $\lam$
is the (discrete fibration version of the) natural transformation whose component at $x\in X$ is
\[
\comp(x/p) \too X(x,x_0) \qv [x\tto^h pa] \mapstoo \lam_a\circ h 
\]
\item
In general if $\lam:m\to\ne(x_0)$ is a cone with a discrete fibration $m:M\to X$ as its base,
the direct image $\ex_f\lam : \ex_f m\to\ex_f\ne(x_0)\iso\ne(fx_0)$ is the natural transformation
whose component at $y\in Y$ is
\[
\comp(y/f\times m) \too \comp(y/f\times X/x_0) \iso Y(y,fx_0) 
\]
\[
[y\tto^h fx, a\in mx] \mapstoo [y\tto^h fx, x=ma\tto^{\lam_a}x_0] \mapstto^\sim f\lam_a\circ h
\]
\item
If $m = \,\down p$, under the correspondence
\[
\ex_f\down p \iso \,\down(f_!p) \qv [y\tto^h fx, [x\tto^k pa]] \mapstto^\sim [y\tto^{fk\circ h}fpa]
\]
the image of the kernel cone corresponds to
\[
[y\tto^h fpa] \mapstoo f\lam_a\circ h	
\]
\item
Since composing the above cone with the unit
\[
f_!p \too \,\down(f_!p) \qv a\in P \mapstoo [fpa\tto^\id fpa]
\]
gives exactly the \v classical" image cone $a\mapstoo f\lam_a$,
the proof is complete.
\end{itemize}
\epf
\noindent So, while the \v classical" image cone has not a straightforward correspective
in the general context, the above proposition suggests to define it as the composite
\[ f_! p \,\too\,\, \down(f_! p) \iso \ex_f\down p \,\too\, \ne(fx_0) \]
\begin{corol}
In $\Cat$, the colimit preserving maps are the classical colimit preseving functors.
\end{corol}
\epf

\subsection{Absolute colimits}

It may happen that the discrete reflection of a space $p:P\to X$ over $X$ gives already a neighborhood:
\eq   \label{13gen} 
\down p \iso \ne(x)
\eeq
so that the reflection map $p\to\,\down p$ is a colimiting cone.
In this case, we say that $x$ is an {\bf absolute} colimit of $p$.
(In particular, the inclusion $x\to \ne(x)$ of a point in its neighborhood is the absolute colimit of $x:1\to X$).
In $\Cat$, the above characterization of absolute colimits was among the leading motifs in~\cite{par73}.

\subsection{Universal displacements}

Given a map $f:X\to Y$ and a point $y:1\to Y$ of its codomain,
a {\bf universal displacement} from $f$ to $y$ is a final point $e:1\to\De_f\ne(y)$  
of the pullback of the neighborhood of $y$: 
\eq   \label{14gen}  
\xymatrix@R=4pc@C=4pc{
1  \ar[r]^e\ar[dr]_x       & M \ar[r]^l\ar[d]^{\De_f\ne(y)} & N \ar[d]^{\ne(y)} & 1 \ar[l]_{e'}\ar[dl]^y \\
                           & X \ar[r]^f                     & Y                 &                         }
\eeq
(in particular, it gives a displacement $l\circ e$ from $fx$ to $y$).

In $\Cat$, $\De_f\ne(y) = f^*\down y $ is the \v comma" or map category $(f\down\,y)$ or $f/y$, 
and we get the usual notion of universal arrow.

\subsection{Adjoint maps}

Given the maps $f:X\to Y$ and $g:Y\to X$ in $\C$, we say that $g$ is a {\bf right adjoint} of $f$ when
\eq    \label{15gen} 
\De_f\adj\De_g:\MX\to\MY 
\eeq
As in~\cite{pt}, we say that $f$ is (right) {\bf adjunctible} 
if for any point $y:1\to Y$ there is a universal displacement
from $f$ to $y$; or, equivalently, if $\De_f$ preserves the neighborhoods (up to isomorphisms), 
that is it restricts to a functor 
\eq    \label{16gen} 
\de_f:\ov Y\to\ov X 
\eeq
with $\ov f \adj \de_f$.
Note that if $f$ has a right adjoint then $\De_f \iso \ex_g$, and so by~(\ref{5gen}) $f$ is also adjunctible.
It is a classical fact that in $\Cat$ the two notions coincide.

\subsection{Dense maps}

We say that a map $f:X\to Y$ in $\C$ is {\bf dense} at the point $y\in Y$ 
if the counit
\eq    \label{17gen} 
\ex_f\De_f\ne(y) \to \ne(y)
\eeq
is a colimiting cone.
(Or, equivalently, if the counit $f_!f^*\ne(y) \to \ne(y)$ is a colimiting cone.)
If $f$ is dense at any point of the codomain we say that it is dense.
Equivalently, $f$ is dense if the composite of the adjoint functors $\,\colim\circ\ex_f \adj \De_f\circ j_X$
\eq   \label{18gen} 
\xymatrix@1@C=3.5pc{\ov Y \ar[r]^j & \MY \ar[r]^{\De_f} & \MX \ar[r]^{\ex_f} & \MY \ar[r]^{\colim} & \ov Y } 
\eeq
is isomorphic to the identity on $\ov Y$; and $f$ is dense at $y$ if
\eq   \label{18cgen} 
\colim(\ex_f\De_f\ne(y)) \iso y 
\eeq
Or, by~(\ref{10bgen}), also
\eq   \label{18bgen} 
\Colim(f\circ\De_f\ne(y)) \iso y
\eeq
Thus, in $\Cat$ we get the usual expression of an object $y\in Y$ as the colimit of
\[
\xymatrix@1@C=3.5pc{f/y \ar[r] & X \ar[r]^f & Y } 
\]
Note also that, by~(\ref{18cgen}), if $f$ is adjunctible then
it is dense iff $\de_f:\ov Y\to\ov X$ is a right inverse of $\ov f:\ov X\to\ov Y$.

\subsection{Fully-faithful maps}

We say that a map $f:X\to Y$ in $\C$ is {\bf fully-faithful} at the point $x\in X$ if the unit 
\eq    \label{17bgen} 
\ne(x) \to \De_f\ex_f\ne(x) = \De_f\ne(fx) 
\eeq
is an isomorphism.
If $f$ is fully-faithful at any point of the domain we say that it is fully-faithful.
Equivalently, $f$ is fully-faithful if the diagram below commutes
in $\Cat$ (up to isomorphisms): 
\eq     \label{19gen} 
\xymatrix@R=4pc@C=4pc{
\ov X  \ar[d]_j \ar[r]^{\ov f} & \ov Y \ar[d]^j    \\
\MX                            & \MY\ar[l]_{\De_f}   }
\eeq
In $\Cat$,~(\ref{17bgen}) expresses, in terms of presheaves, the natural transformation
\[ X(-,x) \to Y(f-,fx) \]
so that we get the usual full and faithful functors.

\section{Some theorems}
\label{thm}

We now show that the above generalized notions are adequate: some classical categorical
theorems can be proved in the general context of $\EM$-categories, finding therein a very 
natural setting.

Some of them become almost tautological, as the following:
\begin{prop}   \label{thm0} 
If the space $X$ has a final point $x$, then a map $p:P\to X$ is final iff $x$
is the absolute colimit of $p$.
\end{prop}
\epf
Others follow easily from the fact that the colimit of a space $p:P\to X$ depends only on 
its discrete reflection, and that if the discrete reflection of $p$ and $q$ are isomorphic, 
then the same remains true after composing them with any map $f:X\to Y$
(see~(\ref{3fs}),~(\ref{5fs}),~(\ref{6fs}) and~(\ref{10bgen})). 
We give the proof of the first one only, the others being equally simple:
\begin{prop} \label{thm1} 
If $x$ is the absolute colimit of $p:P\to X$, then for any $f:X\to Y$, $fx$
is the absolute colimit of $fp:P\to Y$.
\end{prop}
\pf
\[
\down(f_!p) \iso \ex_f \down p \iso \ex_f \down x \iso\, \down(f_!x)
\]
\epf
\begin{prop}   \label{thm2} 
If two spaces $p:P\to X$ and $q:Q\to X$ over $X$ have isomorphic reflections 
$\,\down p\,\iso\,\down q\,$, then for any $f:X\to Y$
\[ \Colim(f\circ p) \iso \Colim(f\circ q) \]
either side existing if the other one does.
\end{prop}
\epf
\begin{corol}   \label{thm3} 
If $e:P\to X$ is a final map, then for any $f:X\to Y$
\[ 
\Colim\,f \iso \Colim(f\circ e)
\]
either side existing if the other does.
\end{corol}
\epf
\noindent 
In Section~\ref{pow} we show that the converse of these propositions 
holds if $X$ has a power object in $\C$.
\begin{corol}    \label{thm4} 
If $P$ has a final point $e:1\to P$, then for any $p:P\to X$
\[ 
\Colim\,p \iso\, pe
\]
and the colimit is absolute.
\end{corol}
\epf

The next two propositions depend both on the following general lemma:
\begin{lemma}   \label{thm5} 
If an object $X$ of a category has a reflection $u:X\to X'$ in a full subcategory,
and if $u$ has a retraction, then the latter is in fact the inverse of $u$. 
\end{lemma} 
\pf
$u\circ r = \id_{X'}$, since both of them solve the universal problem of factoring $u$ through itself:
\[ 
\xymatrix@1@C=3.5pc{X\ar[r]^u\ar@/_1.2pc/[rr]_\id & X' \ar[r]^r\ar@/^1.2pc/@{..>}[rr]^\id & X\ar[r]^u & X' } 
\]
\epf

\begin{prop}    \label{thm6} 
If a final map $e:P\to X$ has a colimit, then it is absolute and is a final point.
\end{prop}
\pf
If $e = \id_X$, it is terminal in $\MX$, and so the colimit $1_X\to \ne(x)$ has a retraction.
Then the lemma applies and $\ne(x) \iso 1_X \iso \,\down \,1_X$, that is $x$ is a final point and
it is the absolute colimit of the identity. 
By Corollary~\ref{thm3}, the same is true for any final map $e$. 
\epf
\begin{corol}    \label{thm7} 
A $\C$-space $X$ has a final point iff some final map has a colimit in it.
\end{corol}
\epf
\begin{prop}[Formal criterion for the existence of universal arrows]   \label{thm8} 
Given a map $f:X\to Y$ in $\C$ and a point $y:1\to Y$, there is a universal displacement from $f$ to $y$ 
iff the discrete space $\De_f \ne(y)$ over $X$ has a colimit, preserved by $f$ itself.
\end{prop}
\pf
One direction is trivial, since any point is the (absolute) colimit of its neighborhood.
Suppose conversely that
\[ 
\lam: \De_f \ne(y) \to \ne(x)
\]
is a colimiting cone. By the lemma, the proposition is proved if we find a retraction of $\lam$.
We claim that such a retraction is given by the adjunct 
\[
u^\star: \ne(x) \to \De_f \ne(y)
\]
of the map $u$, universally induced by the image cone, which is itself colimiting by hypotesis: 
\eq  
\xymatrix@R=3pc@C=4pc{
\ex_f \De_f \ne(y) \ar[d]^{\ex_f\lam} \ar[dr]^\eps   &     \\
\ex_f \ne(x) \ar@{..>}[r]^-u                         &   \ne(y)   \\   }                                   
\eeq
(where $\eps$ is the counit of the adjunction $\De_f \adj \ex_f$).
Indeed by the naturality of the adjunction bijection
\[
(-)^\star : \MY(\ex_f \_\_\, , -) \to \MX( \_\_\, , \De_f - ) \qq 
\]
we have 
\[
u^\star\circ\lam = (u\circ\ex_f\lam)^\star = \eps^\star = \id_{\De_f \ne(y)}
\] 
\epf

\begin{prop}   \label{thm9} 
Any adjunctible map is colimit preserving.
\end{prop}
\pf
Consider $f:X\to Y$ such that $\De_f$ preserves neighborhoods, and a colimiting 
cone $\lam:m\to \ne(x)$ over $X$ with a discrete base.
We want to show that $\ex_f\lam$ is itself colimiting: 
for any $l:\ex_f m\to \ne(x')$ over $Y$ there is a unique $v$ in the left hand diagram below: 
\[
\xymatrix@R=2.5pc@C=4pc{
\ex_f m \ar[r]^{\ex_f\lam}\ar[dr]^l & \ex_f \ne(x) \ar@{..>}[d]^v  \\
                                    & \ne(x')                         }
\qq\qq
\xymatrix@R=2.5pc@C=4pc{
m \ar[r]^\lam\ar[dr]^{l^\star} & \ne(x) \ar@{..>}[d]^u  \\
                               & \De_f \ne(x')                 }
\]
By hypotesis, $\De_f \ne(x')$ is itself a neighborhood, so that we have a universally
induced $u$ in the right hand diagram.
Now, using again the naturality of the adjunction bijections
\[
(-)_\star : \MX( \_\_\, , \De_f - ) \to \MY(\ex_f \_\_\, , -) \qv (-)^\star = (-)_\star^{-1}
\]
one easily checks the unicity:
\[
l = v\circ\ex_f\lam \imp l^\star = (v\circ\ex_f\lam)^\star = v^\star\lam \imp v^\star = u \imp v = u_\star 
\]
and the existence:
\[
u_\star\circ\ex_f\lam = (u\circ\lam)_\star = (l^\star)_\star = l
\]
\epf

\section{Power objects}
\label{pow}

So far, we have only used the $\EM$-structure of $\Cat$.
Now we tentatively suggest how other aspects of the rich structure of $\Cat$ may be exploited, 
by introducing further axioms on $\C$.
We begin with power objects, playing the role of the presheaf categories in $\Cat$.

We say that
\[ \y:X\to\P X \]
is a {\bf Yoneda map} if the following adjoint composites constitute 
an equivalence of categories $\MX \equ \ov{\P X}$ between the category of discrete spaces over $X$ and
the adherence category of $\P X$:
\[
\xymatrix@1@C=3.5pc{\MX \ar[r]^-{\ex_\y} & \MPX \ar[r]^-{\colim} & \ov{\P X} } 
\]
\[
\xymatrix@1@C=3.5pc{\ov{\P X} \ar[r]^-j & \MPX \ar[r]^-{\De_\y} & \MX } 
\]
If this is the case, we also say that the codomain $\P X$ of $\y$ is a power object for $X$ (via $\y$). 
In particular $\y: X \to \P X$ is a dense map by definition, and it is also fully-faithful 
as one sees by composing the following commutative diagram (recall~(\ref{5gen}))   
\eq  
\xymatrix@R=2pc@C=4pc{
\ov X  \ar[d]_j \ar@/^1pc/[drr]^{\ov\y} &     \\
\MX   \ar[r]^{\ex_\y} & \MPX \ar[r]^\colim  & \ov{\P X}   }
\eeq
with $\xymatrix@1@C=3.5pc{\ov{\P X} \ar[r]^-j & \MPX \ar[r]^-{\De_\y} & \MX }$. 

The fact that in $\Cat$ the Yoneda embedding $\y:X\to\PrX$ is a Yoneda map as just defined, 
expresses in particular the following well-known facts:
\begin{itemize}
\item
The map $\y/- :\PrX \to \MX$ is an equivelence. This fact, in turn, can be seen as a strong form
of the Yoneda Lemma itself.
\item
Any presheaf is a canonical colimit of representable presheaves (density of $\y$; see~(\ref{18bgen})). 
\end{itemize}
Now we show that if the space $X$ has a power object then the discrete reflection of a space over $X$
can be expressed as a colimit in $\P X$:
\begin{prop}
Let $\y:X\to\P X$ be a Yoneda map. For any $p\in\CX$
\[ 
\down p \iso \De_\y\,j\,\Colim(\y\circ p)
\]
\end{prop}
\pf
Just recall~(\ref{10bgen}):
\[
\down p \iso (\De_\y\circ j\circ\colim\circ\ex_\y)\down p \iso \De_\y\,j\,\colim(\ex_\y\down p) \iso \De_\y\,j\,\Colim(\y\circ p)
\]
\epf
The following proposition was proved in~\cite{par73} in the case $\C = \Cat$.
\begin{corol} 
If the space $X$ has a power object, 
two spaces $p:P\to X$ and $q:Q\to X$ over $X$ have isomorphic discrete reflections 
$\,\down p\,\iso\,\down q\,$, iff for any $f:X\to Y$
\[ \Colim(f\circ p) \iso \Colim(f\circ q) \]
either side existing if the other one does. In particular, 
\begin{itemize}
\item
a map $e:P\to X$ is final iff for any $f:X\to Y$
\[ 
\Colim(f\circ e) \iso \Colim\, f  
\]
either side existing if the other does.
\item
the point $x$ is the absolute colimit of $p:P\to X$ iff for any $f:X\to Y$
\[ 
\Colim (f\circ p) \iso fx
\]
\end{itemize}
\end{corol}
\pf
One direction was proved in Proposition~\ref{thm2}.
The other one follows from the above proposition with $f = \y$.
For the particular cases, consider $q = \id_X$ and $q = x$ respectively.  
\epf

\section{Duality and exponentials}
\label{dual}

What we have done so far, could be called \v one-sided category theory" (say, \v left-sided"), 
since it is modeled on one of the two comprehensive factorization systems:
final functors and discrete fibrations, rather than initial functors and discrete opfibrations.
So we have defined colimits rather than limits.
But also in this context we can define the product of two points as the following
universal displacement:
\[ 
\xymatrix@R=4pc@C=4pc{
1  \ar[r]^e\ar[dr]_{x\times y} & M \ar[r]^l\ar[d]^{f^*n} & N \ar[d]^n & 1 \ar[l]_{e'}\ar[dl]^{\la x , y \ra} \\
                               & X \ar[r]^\De            & X\times X  &                      }
\]
Similarly, if $\C$ is cartesian closed we can define the limit of any map.

However, to have a balanced theory we need to assume a duality on $\C$,
modeled on the functor $(-)\op:\Cat\to\Cat$.

We say that a $\EM$-category $\C$ is \v two-sided" if it is equipped with an isomorphism
\[
(-)':\C\to\C
\]
Then any concept and property so far presented can be dualized.
For example, $e:X\to Y$ is an \v initial map" iff $e':X'\to Y'$ is final;
similarly, one says that the space $p:P\to X$ is \v right discrete" over $X$ whenever $p':P'\to X'$
is (left) discrete.
So we have another factorization system $(\E',\M')$ with the corresponding functors
\[ 
\up\,\,\,\adj i':\MoX\to\CX
\]
and for any $f:X\to Y$,
\[
\ex_f'\adj\De_f':\MoY\to\MoX
\]
There are the \v right" adherence categories $\ov X'$ and the partially defined \v limit" functors:
\[
\lim\adj j':\ov X' \to \MoX \qv \Lim\adj k':\ov X' \to \CX
\]
The functor $(-)':\C\to\C$ induces isomorphisms
\[
(-)':\CX\to\CX'
\]
which restrict to isomorphisms
\[
\MX\iso\MoX'
\] 
or equivalently
\[
\MX'\iso\MoX
\] 
and also
\[
\ov X'\iso\ov{X'}
\]
Under these equivalences, $\De_{f'}$ corresponds to $\De_f'$, $\ex_{f'}$ to $\ex_f'$, $\colim_{X'}$ to $\lim_X$,
and so on.

\subsection{Axioms on duality}

It seems appropriate to require that the duality functor fixes 
(at least up to isomorphisms) the subcategory $\S\inc\C$ of (left) sets.
In particular, left and right sets coincide, and the constant spaces over a space 
are left and right discrete (see~(\ref{4gen})). 
If $\C$ is two-sided, it seems also natural to require that 
the adherence functor preserves duality: in particular,
\[ \ov{X'} \iso \ov X \op \]

\subsection{Hom maps}

If there is a Yoneda map $\omega:1\to\Omega$,
we say that $\Omega$ is a {\bf truth values} or a {\bf internal sets} object.
In particular $\ov\Omega \equ \S$.

Suppose that $\C$ is cartesian closed and has a truth values object.
A {\bf hom map} for $X\in\C$ is a map 
\[ X'\times X \to \Omega \]
such that both the transposes
\[ X\to \Omega^{X'} \qv X'\to \Omega^X \]
are Yoneda maps.
When $X$ has a hom map, we also say that $X$ is {\bf locally small}.

\subsection{Axioms on exponentials and the reflection formula}
\label{exp}

Although the category of categories is (cartesian closed but) not locally cartesian closed,
some interesting laws concerning exponentials in $\CatX$ hold (see~\cite{pis}); 
we now assume that they hold in the two-sided $\EM$-category $\C$:
\begin{itemize}
\item
Left or right discrete spaces over $X$ are exponentiable in $\CX$.
\item
If $m\in\MX$ and $n\in\MoX$ are left and right discrete spaces over $X$,
then $m^n\in\MX$ and, symmetrically, $n^m\in\MoX$.
\end{itemize}
In particular, since constant maps are left and right discrete, 
we can define a \v complement operator" parameterized by the sets $S\in\S$ 
\[ \neg : (\MX)\op \to \MoX \qv m\mapsto (\pi_S)^m \]
where the projection $\pi_S : X\times S \to X$ is the constant space $(!_X)^*S$ over $X$ with value $S$;
and we can prove the following law as in~\cite{pis}: 
\begin{prop}  \label{dual1}
For any $p,q\in\CX$,
\( \comp(p^*\down q) \iso \comp(q^*\up p) \)
\end{prop}
\pf
Here, $\comp := \comp\circ\tot : \MX\to\S$ gives the components $\comp M = \comp(\tot\, m)$
of the total of $m:M\to X$. So, its right adjoint takes $S\in\S$ into $\pi_S$.
\[
\begin{array}{c}
\comp(p^*\down q) \to S \\ \hline 
\comp(p\,\times\!\down q) \to S \\ \hline 
p\,\times\!\down q \to\,\, \pi_S \\ \hline 
p \to\,(\pi_S)^{\down\,\,q} \\ \hline 
\up p\to\,(\pi_S)^{\down\,\,q} \\ \hline 
\comp(\up p\,\times\!\down q) \to S \\ \hline 
\comp(\up p\times q) \to S 
\end{array}
\]
\epf
In particular, we get the \v reflection formula", which gives the value
(see~(\ref{4gen})) of the discrete reflection $\down p$ at any point of $x$ of $X$:
\eq  \label{1dual}
(\down p)x \iso \comp(p^*\up x) 
\eeq

\section{Categories of categories}
\label{ar}

Although we have developed a certain amount of category theory in the general context 
of a $\EM$-category $\C$, possibly enriched with power objects or a duality functor, 
we have already stressed that in general the objects of $\C$ should be considered more basic 
than categories.
It is the presence of an arrow object, playing the role of the arrow category $2\in\Cat$,
which makes $\C$ more concretely a \v category of categories", allowing an analysis
\v inside" the objects $X\in\C$, through their \v arrows" $2\to X$.

In the early sixties, Lawvere showed how the arrow category can be taken as a base 
for axiomatizing the category of categories.
Assuming the finite completeness and cocompleteness of $\C$, along with 
suitable axioms on $2$ and on the pushouts $3$ and $4$, the geometric way of analyzing 
the objects of $\C$ via the figures with these shapes shows that any object $X\in\C$ \v is" a category. 
In~\cite{law63} and~\cite{law66}, sets $S\in\C$ are defined as those objects orthogonal to $2\to 1$
(that is, any map $2\to S$ is constant), and it is assumed that sets are reflective and coreflective in $\C$.
(It is also assumed that $\C$ is cartesian closed, so that any $X\in\C$ is in fact a $\S$-category.)  

In our context, discrete objects (over any object) are given as part of the structure, 
and the reflection therein (over any object) is the key axiom. 
If $\C$ is a two-sided $\EM$-category, an arrow object is defined in terms of discreteness,
rather than the other way round:
it is a bipointed object of $\C$ which generates the factorization system and its dual.

Here, we explore some convergences between the two approaches.
But before considering the bipointed object $2$, we begin by assuming a pointed object $t:1\to T$
which generates the factorization system.
If $\C$ is a \v topological category" (see Section~\ref{top}),  
such an object should rather be thought of as a \v universal convergence object".
It may play a role similar to the \v infinitesimal" pointed object of Lawvere's cohesive categories.
Indeed, the analogy between such a $T\in\C$ and the arrow category $2\in\Cat$ has been
pointed out in~\cite{law03}.

\subsection{The universal convergence object}
\label{conv}

We now assume that $\EM$ is generated by a pointed object $t:1\to T$, that is $\M = t^\perp$.
In particular $t$ is final, so that $T$ is connected (that is, $T\to 1$ is itself final). 
First, we show that $S\in\C$ is a set (in our sense) iff it is a set in the sense 
of~\cite{law66} and~\cite{law03}:
\begin{prop}
A $\C$-space $S$ is a set iff any $T\to S$ is constant.
\end{prop}
\pf
In one direction, if $S$ is a set then $T\to 1$ is orthogonal to $S\to 1$:
\[
\xymatrix@R=3pc@C=3pc{
T \ar[r]^f\ar[d]          &  S \ar[d] \\
1 \ar[r]\ar@{..>}[ur]^u   &  1         }
\]
so that any map $f:T\to S$ is constant.
In the other direction, suppose that the latter is true.
We want to show that $t:1\to T$ is orthogonal to $S\to 1$:
\[
\xymatrix@R=3pc@C=3pc{
1 \ar[r]^x\ar[d]_t        &  S \ar[d] \\
T \ar[r]\ar@{..>}[ur]^u   &  1         }
\]
Unicity: since $u$ is constant, we have $u = y\circ !$; and since the upper triangle 
commutes, $x = u\circ t = y\circ !\circ t = y$, so that $u = x\circ !$.
Existence: $x\circ !\circ t = x$, while also the lower triangle clearly commutes.
\epf

In a \v topological category", a map $f:T\to X$ from the \v universal convergence 
object" $t:1\to T$ is \v a convergence" to the point $ft$.
Thus, a space $m:M\to X$ over $X$ is discrete iff
any convergence to a point of the base $X$ has a unique lifting to a convergence 
to any point over it:
\[
\xymatrix@R=3pc@C=3pc{
1 \ar[r]^a\ar[d]_t          &  M \ar[d]^m \\
T \ar[r]^f\ar@{..>}[ur]^u   &  X           }
\]
(In particular, if $P$ has no points, any space $p:P\to X$ over $X$ is discrete.)
More specifically, referring to the \v classical" analysis of infinitesimal aspects in topology 
via ultrafilter convergence, $t:1\to T$ may be seen as the \v free converging ultrafilter". 
Then a map $f:T\to X$ becomes an \v ultrafilter" converging to $ft$ and we get
a characterization of discrete spaces over $X$ as \v discrete ultrafilter fibrations" 
in the spirit of~\cite{hof05}.

Note that, by Corollary~\ref{thm4}, if $f:T\to X$ is a convergence to $x = ft$ in the above sense 
then $x$ is an absolute colimit of $f$.
In fact we can be more precise.
Since in the square below $e$ and $t$ are both final points, the unique diagonal is also final,
so that the neighborhood of $x$ is the discrete reflection of $f$, and $u$ is the (absolute)
colimiting cone:
\[
\xymatrix@R=3pc@C=3pc{
1 \ar[r]^e\ar[d]_t         &  N \ar[d]^{\ne(x)} \\
T \ar[r]^f\ar@{..>}[ur]^u  &  X                    }
\]

\subsection{The arrow object}

If the convergence object is in fact bipointed, we denote it by $s,t:1\to 2$
and call it a {\bf left arrow object} for $\C$.
If $\EM$ is two-sided and $s$ generates the dual of $\EM$, $2$ is an {\bf arrow object}.
A map $l:2\to X$ is an {\bf arrow} of $X$, \v from the (domain) point $x = l\circ s$ 
to the (codomain) point $y = l\circ t$". We also write $l : x \imp y$. 
Given a point $x:1\to X$, the constant arrow $l = x\circ !$ at $x$ is also called
the \v identity" $l = \id_x : x \imp x$ of $x$.

Thus, if $\C$ has a (left) arrow object, a space $m:M\to X$ over $X$ is discrete iff
any arrow of the base $X$ has a unique lifting to a convergence 
to any point over its codomain:
\[
\xymatrix@R=3pc@C=3pc{
1 \ar[r]^a\ar[d]_t          &  M \ar[d]^m \\
2 \ar[r]^l\ar@{..>}[ur]^u   &  X           }
\]
Since, as before, $2\to 1$ is also final (that is, $2$ is connected) we also have 
\[
\xymatrix@R=3pc@C=3pc{
2 \ar[r]^l\ar[d]            &  M \ar[d]^m \\
1 \ar[r]^x\ar@{..>}[ur]^u   &  X           }
\]
That is, any arrow which is projected to an identity by a discrete space
is itself an identity.

Suppose now that Lawvere's axioms on $2$, $3$ and $4$ hold.
In particular, the pushout $3$ 
\[
\xymatrix@R=3pc@C=3pc{
1 \ar[r]^s\ar[d]_t     &  2 \ar[d]^{l_2} \\
2 \ar[r]^{l_1}         &  3              }
\]
has three non-identity arrows, and its \v internal picture"
is the following:
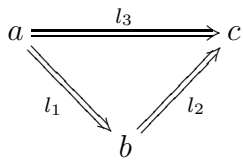
\[
\xymatrix@R=2.5pc@C=2.5pc{
a \ar@{=>}[dr]_{l_1}\ar@{=>}[rr]^{l_3} &                        &  c    \\
                                       & b \ar@{=>}[ur]_{l_2}   &        }
\] 
Given two consecutive arrows $x\imp y\imp z$ of $X$, their internal composite is then the
composite $u\circ l_3$, where $u:3\to X$ is the uniquely induced map from the pushout. 

In this way, to any $\C$-space $X\in\C$ there corresponds a category $X^\star$, 
(whose identities are the constant arrows as above), and to any map $f:X\to Y$ in $\C$
there corresponds a functor $f^\star:X^\star\to Y^\star$ (see~\cite{law66}).
Furthermore, $m:M\to X$ is a discrete space over $X$ iff $m^\star:M^\star\to X^\star$
is a discrete fibration over $X$: the diagram
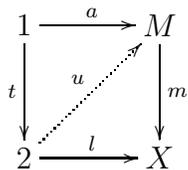
\[
\xymatrix@R=3pc@C=3pc{
1 \ar[r]^a\ar[d]_t          &  M \ar[d]^m \\
2 \ar[r]^l\ar@{..>}[ur]^u   &  X           }
\]
says that any arrow in $X^\star$ has a unique lifting in $M^\star$ to any
object of $\M^\star$ over its codomain.

Thus we have a (discreteness preserving) pseudo-functor $(-)^\star:\C\to\Cat$, and 
we presently show that there is a natural transformation $X^\star\to\ov X:\C\to\Cat$
to the adherence pseudo-functor. 
\begin{lemma}
Let $m:M\to X$ be a discrete space over $X$, and $l'':x\imp z$ be the composite of
two consecutive arrows $l:x\imp y$ and $l':y\imp z$ of $X$.
Then the domain of the lifting of $l''$ to a point $a$ over $z$ coincides with
the domain of the lifting of $l$ to the domain of the lifting of $l'$ to $a$. 
\end{lemma}
\pf
First, recall that by Proposition~\ref{fs1} the arrow $l_2:2\to 3$ is final, 
so that also $l_2\circ t : 1 \to 3$ is a final point.
Then, given a discrete space $m:M\to X$ over $X$, and two consecutive arrows $l:x\imp y$
and $l':y\imp z$ of $X$, the lifting of their composite $l'':x\imp z$ to a point over $z$ 
can be obtained by composing $l_3$ with the lifting of the corresponding \v triangle" $3\to X$:
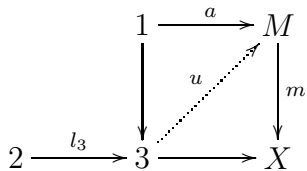
\[
\xymatrix@R=3pc@C=3pc{
              & 1 \ar[r]^a\ar[d]          &  M \ar[d]^m \\
2\ar[r]^{l_3} & 3 \ar[r]\ar@{..>}[ur]^u   &  X           }
\]
The lemma now follows.
\epf
Now, to any arrow $l:x\imp y$ in $X$ we associate a displacement $\lam$ from $x$ to $y$
(see Section~\ref{disp}), that is an arrow of the adherence category $\ov X$.
This is given by the domain $u\circ s$ of the lifting of $l$ to the final point of (the total of)
the neighborhood $\ne(y)$:
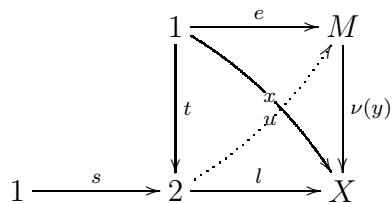
\[
\xymatrix@R=4pc@C=4pc{
          & 1  \ar[d]^t\ar[r]^e\ar@/^/[dr]|x   & M \ar[d]^{\ne(y)} \\
1\ar[r]^s & 2  \ar[r]^l\ar@{.>}@/_/[ur]|u      & X                  }
\]
\begin{prop}
The above assignment defines a functor $X^\star\to\ov X$.
\end{prop}
\pf
Since identities are easily seen to be preserved, we concentrate on composition.
Let $l:x\imp y$ and $l':y\imp z$ in $X^\star$, and let $\lam$ and $\lam'$ be the corresponding
arrows in $\ov X$:
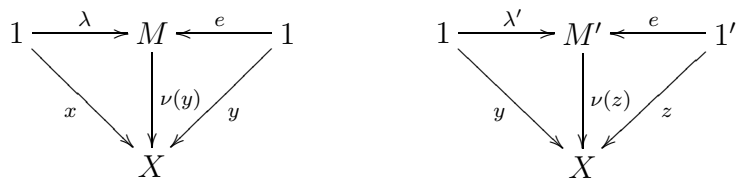
\[
\xymatrix@R=3pc@C=3pc{
1  \ar[r]^\lam\ar[dr]_x    & M \ar[d]^{\ne(y)} & 1 \ar[l]_e\ar[dl]^y \\
                           & X                 &                        }
\qq\qq
\xymatrix@R=3pc@C=3pc{
1  \ar[r]^{\lam'}\ar[dr]_y    & M' \ar[d]^{\ne(z)} & 1 \ar[l]_e'\ar[dl]^z \\
                              & X                  &                        }
\]
Their composite in $\ov X$ is the composite
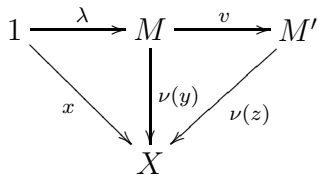
\[
\xymatrix@R=3pc@C=3pc{
1  \ar[r]^\lam\ar[dr]_x    & M \ar[r]^v\ar[d]^{\ne(y)} & M' \ar[dl]^{\ne(z)} \\
                           & X                         &                        }
\]
where $v$ is universally induced by $\lam'$.
But if $h:2\to\ne(y)$ is the lifting of $l:2\to X$ to $e$, then $v\circ h:2\to\ne(z)$ 
is the lifting of $l$ to $v\circ e = h'\circ s$, the domain of the lifting of $l'$.
Thus, the proposition follows from the lemma above.
\epf
\begin{prop}
The functors $X^\star\to\ov X$ are the components of a pseudo-natural transformation.
\end{prop}
\pf
Let $f:X\to Y$ be a map in $\C$. Suppose that $l:x\imp x'$ is an arrow in $X$ 
and that $fx = y$ and $fx' = y'$.
We must show that if the displacement $\lam:\ne(x)\to\ne(x')$ corresponds to the 
arrow $l:2\to X$, then the displacement $\ov f\lam = \ex_f\lam : \ex_f\ne(x)\to\ex_f\ne(x')\iso\ne(y')$ 
corresponds to the arrow $f l:2\to X$.
Comparing the following diagrams:
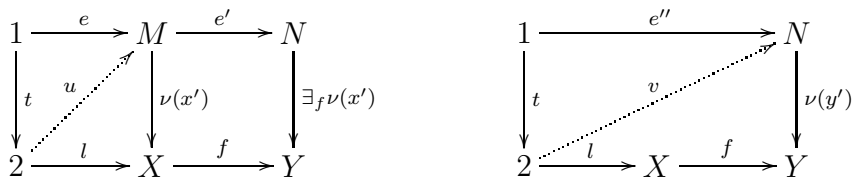
\[
\xymatrix@R=3pc@C=3pc{
1\ar[r]^e\ar[d]^t           & M \ar[r]^{e'}\ar[d]^{\ne(x')}   &  N \ar[d]^{\ex_f\ne(x')} \\
2\ar[r]^l\ar@{..>}[ur]^u    & X \ar[r]^f                      &  Y                        }
\qq\qq
\xymatrix@R=3pc@C=3pc{
1\ar[rr]^{e''}\ar[d]^t      &               &  N \ar[d]^{\ne(y')} \\
2\ar[r]^l\ar@{..>}[urr]^v   & X \ar[r]^f    &  Y                   }
\]
one sees that $v = e'\circ u$, and the proposition follows.
\epf

\section{Further examples}
\label{ex}

\subsection{\v Classical" factorization systems}

As hinted in the Introduction, if the points are discrete (as it happens for most
of the \v classical" factorization systems) then the concepts defined in
Section~\ref{gen} (that is, $\EM$-category theory) tend to became trivial.

For instance, let $\EM$ be the classical (strong) epi-mono factorization system on $\Set$
or on another category $\C$.
Then, for any $X\in\C$, $\ov X$ is discrete. 
Thus, a map $p:P\to X$ has a colimit iff it is strongly constant, that is factors {\em uniquely} 
through $1$. (If $X$ has more than a point, $0\to X$ has not a colimit.) 
If $X$ has a subobject different from the maximum $\id_X:X\to X$, 
then $\MX$ is not discrete, so that $X$ cannot have a power object.
If $\C = \Set$, a map is adjunctible iff it is bijective.
A map is dense iff it is surjective, and fully-faithful iff it is injective.

Similar results hold for factorization systems where all the arrows are in $\M$ 
(so that $\E$ is the class of isomorphisms).

\subsection{Graphs}

Let $\Gph$ be the category of reflexive graphs and $s,t:1\to 2$ the arrow graph.
The pre-factorization system $\EM$ generated by $t$ is in fact a 
factorization system (see~\cite{pis} and references therein).
So, $\Gph$ is a two-sided $\EM$-category with an arrow object.
The category $\Gph$ is \v incomplete" in two related respects:
\begin{enumerate}
\item
there is no internal composition of internal arrows, so that $\ov X$ is a sort of \v completion"
of the graph $X$: indeed, it is the free category on it;
\item
not all functors $\ov X \to \ov Y$ are realized by an actual map $f:X\to Y$ in $\C$.
\end{enumerate}
In $\Gph$ there are no power objects. Indeed, if $\y:X\to\P X$ would exist, then $\MX \iso \ov{\P X}$; 
but the latter should be a free category, while $\MX \iso \Set^{\ov X\op}$ is not so.
Colimits in a graph $X$ are the same that colimits in the corresponding free category.

Similarly, if $\Gph'$ is the category of {\em irreflexive} graphs
and $s,t:1\to 2$ the arrow irreflexive graph,
then the pre-factorization system $\EM$ generated by $t$ is a 
factorization system. Note however that:
\begin{itemize}
\item
the arrow irreflexive graph $2\in\Gph'$ is not an arrow object (it has no points);
\item
the $\Gph'$-sets are the sets with an endomapping;
\item
the connected $\Gph'$-spaces are those graphs whose endomapping
reflection is terminal;
for example, the dot graph and $2$ itself are not connected (see~\cite{pis}).
\end{itemize}

\subsection{Posets}
\label{pos}

Let $\Pos$ be the category of posets, $s,t:1\to 2$ the arrow poset,
and $\EM$ the pre-factorization system generated by $t$ and $1+1\to 1$. 
A map $f:X\to Y$ in $\Pos$ is orthogonal to $1+1\to 1$ iff it is injective on objects,
so that the maps in $\M$ are (isomorphic to) the lower-sets inclusions.
So, $\MX$ is the category (in fact, poset) $\down X$ of lower-sets of $X$, and since
the latter is reflective in $\Pos/X$, $\EM$ is a factorization system.
Furthermore, a map $f:X\to Y$ is in $\E$ iff its lower-set reflection is the maximum $Y$,
that is iff it is a \v cofinal mapping" in the classical sense. 
 
Given  a point $x\in X$ of a poset, its neighborhood $\ne(x) = \,\down x$ is 
the principal lower-set generated by $x$.
As in $\Cat$, $\ov X \iso X$.
The colimit of a map $p:P\to X$ in $\Pos$ is the sup of the 
set of points $pa ,\, a\in P$, or equivalently of the lower-set $\down p$ 
generated by $p$. 
The colimit is absolute iff such a sup is in fact a maximum.
A map $f:X\to Y$ is dense if any $y\in Y$ is the sup of the $fx$ which are
less than $y$ (that is, iff the functor $f$ is dense in the classical sense).

The category of $\Pos$-sets is $\S = 2$, and
the component functor $\Pos\to 2$ reduces to the \v non-void" predicate.
Every $X\in\Pos$ has a power object: the Yoneda map $\y : X\to\,\,\down X$
is the inclusion of $X$ in its lower-sets poset.
In particular, the truth-values poset is $\Omega = 2$,
and we have the natural hom maps $X'\times X \to 2$.
Although $2\in\Pos$ clearly satisfies Lawvere's axioms, 
it is not an arrow object in our sense, since it does not generate $\EM$.

Note that another factorization system on $\Pos$ can be obtained by simply
restricting the comprehensive factorization system to the full subcategory of posets.
(Indeed, if $m:M\to X$ is a discrete fibration and $X$ is a poset, then also $M$ is a poset.)
The neighborhoods are again the principal lower-sets, but now there are no power objects.
This $\EM$-category may be considered as representing the common kernel shared by 
\v categories of categories" and \v topological categories".

\subsection{Topological spaces?}
\label{top}

As already remarked, by exploiting the logic of a factorization system on an
arbitrary finitely complete category, we have been doing in a sense 
$\EM$-category theory and $\EM$-topology at the same time.
Indeed, we have often let our intuition be guided by topological ideas, and we have also 
talked about \v topological categories" when we meant to refer specifically 
to some category of \v spaces" rather than of \v categories".

Of course, we are not thinking of the classical category of topological spaces,
since we assume that some kind of infinitesimal objects exist in $\C$. 
In particular the neighborhood $\ne(x)$ of a point should be thought of as
the infinitesimal part of the space near $x$.

We have also observed in Section~\ref{ar} 
that the classical analysis of infinitesimal aspects
in topology via ultrafilter convergence can be in part reconstructed in our general
context by assuming a \v universal convergence object" $t:1\to T$.
Here, rather than presenting any specific model of \v topological" $\EM$-category,
we briefly speculate on some properties which may be required to hold in it. 

First, it seems natural to assume that neighborhoods are really {\em parts} of the space,
that is, the maps $\ne(x):N\to X$ are monorphisms.
Thus, a space over $X$ has at most a way to be near $x$ (see Section~\ref{disp}).
Although both the $\EM$-categories of posets just defined have this property
they do not have the following \v separation property": distinct points
have disjoint neighborhoods (that is, with an initial pullback). 

If this is the case, there are no non-trivial displacements between points, so that $\ov X$ is discrete.
Thus, $p:P\to X$ (with $P$ non-empty) is near $x$ iff it converges to $x$. 
Furthermore, a map $f:X\to Y$ is dense iff for any $y\in Y$ the neighborhood $\ne(y)$
meets some (infinitesimal) part of $X$.

As observed in Section~\ref{ne}, spaces $T$ with a final point $t:1\to T$
(for instance, the total space of any neighborhood) should be considered as 
\v absolutely concentrated" around their final point $t$.
The neighborhood of $t$ itself is the identity $\ne(t) = \id_T:T\to T$.
The terminal map $T\to 1$ is fully-faithful and adjunctible.
Conversely, if $X\to 1$ is adjunctible, then $X$ has a terminal point.

The monomorphic discrete spaces $m:M\to X$ over $X$ (in particular the neighborhoods) 
can be seen as the \v open parts" of $X$.
Thus any map in $\C$ is continuous in the classical sense. 

One may reasonably assume that $\C$ is a category of cohesion relative to $\S$,
in the sense of~\cite{law07}. We now consider only the existence of a right adjoint
to the discrete inclusion, but {\em at any level}.
That is, we suppose that for any $X\in\C$ there is an \v interior operator" $i\adj (-)^\circ :\CX\to\MX$ 
(as it happens in $\Cat$, but also in the classical category of topological spaces, 
if $\M$ is the class of local homeomorphisms).
In that case, the points of $p^\circ$ over $x\in X$ are bijective with morphisms $\ne(x)\to p$ over $X$.
(This corresponds to the classical fact that the points of the fiber over $x$ of the \'etal\'e  
coreflection of a bundle are bijective with the germs of its cross sections at $x$.) 
In particular, for monomorphic spaces over $X$, a point is in the interior of a part
iff its neighborhood is in the part.

So far, we have done \v one-sided $\EM$-topology".
Clearly, it is not appropriate to assume a duality functor $\C\to\C$ as for \v categories of categories".
Yet, $\C$ may be considered two-sided, by assuming that it has another factorization system $(\E',\M')$,
where $\M'$ plays the role of the class of continuous perfect maps
(so that $\S' = \M'/1$ is the subcategory of \v compact" spaces). 
Indeed, the latter are classical characterized by a \v discrete ultrafilter fibration" property,
dual to the one that holds for local homeomorphisms.
In that case, it seems likely that the law of Proposition~\ref{dual1} may play an
important role.
(As the proof of the proposition shows, the law itself would follow from a complement operator parameterized by $\S$,
which exchanges left and right discrete spaces over $X$; that is assuming, for any $m\in\MX$, a right adjoint to
the functor $\comp(m\times -):\CX\to\S$, valued in $\MoX$, and conversely.) 
Indeed, we would then have in particular the~(\ref{1dual}) and its dual, which,
for monomorphic spaces over $X$, can be rephrased as:
\v a point is in the closure of a part iff its neighborhood meets the part"; 
and as: \v a point is adherent to a part iff its closure meets the part".



\begin{refs}

 


\bibitem[Clementino, Hofmann \& Janelidze 2005]{hof05} M.M. Clementino, D. Hofmann and G. Janelidze (2005), 
Local Homeomorphisms via Ultrafilter Convergence, {\em Proc. Amer. Math. Soc.} {\bf 133}, 917-922.



\bibitem[Kock \& Reyes, 1977]{kor77} A. Kock and G.E. Reyes (1977), {\em Doctrines in Categorical Logic},
Handbook of Mathematical Logic, North-Holland, 283-313.

\bibitem[Lawvere, 1963]{law63} F.W. Lawvere (1963), {\em Functorial Semantics of Algebraic Theories},
Phd Thesis, republished in {\em Reprints in Theory and Appl. Cat.}

\bibitem[Lawvere, 1966]{law66} F.W. Lawvere (1966), {\em The Category of Categories as a Foundation for Mathematics},
Proceedings of the Conference on Categorical Algebra, La Jolla, 1965, Springer, New York, 1-20.

\bibitem[Lawvere, 1969]{law69} F.W. Lawvere (1969), Adjointness in Foundations,
{\em Dialectica} {\bf 23}, 281-295. Republished in {\em Reprints in Theory and Appl. Cat.}
 
\bibitem[Lawvere, 1970]{law70} F.W. Lawvere (1970), {\em Equality in Hyperdoctrines and the Comprehension Scheme as an Adjoint Functor},
Proceedings of the AMS Symposium on Pure Mathematics, XVII, 1-14. 

\bibitem[Lawvere, 2003]{law03} F.W. Lawvere (2003), Foundations and Applications: Axiomatization and Education,
{\em Bull. Symb. Logic} {\bf 9}(2), 213-224.

\bibitem[Lawvere, 2007]{law07} F.W. Lawvere (2007), Axiomatic Cohesion, 
{\em Theory and Appl. Cat.} {\bf 19}, 41-49. 


\bibitem[Mac Lane, 1971]{mac71} S. Mac Lane (1971), {\em Categories for the Working Mathematician}, 
Graduate Texts in Mathematics, vol. 5, Springer, Berlin.


\bibitem[Par\'e, 1973]{par73} R. Par\'e (1973), Connected Components and Colimits,
{\em J. Pure Appl. Algebra} {\bf 3}, 21-42.

\bibitem[Pisani, 2007a]{pis} C. Pisani (2007a), Components, Complements and the Reflection Formula, 
{\em Theory and Appl. Cat.} {\bf 19}, 19-40. 

\bibitem[Pisani, 2007b]{pis2} C. Pisani (2007b), Components, Complements and Reflection Formulas, preprint, math.\-CT/0701457. 

\bibitem[Street \& Walters, 1973]{stw73} R. Street and R.F.C. Walters (1973), The Comprehensive Factorization of a Functor,
{\em Bull. Amer. Math. Soc.} {\bf 79}(2), 936-941.

\bibitem[Taylor, 1999]{pt} P. Taylor (1999), {\em Practical Foundations of Mathematics}, 
Cambridge Studies in Advanced Mathematics, Cambridge University Press.

\end{refs}

\end{document}